\documentclass[12pt]{amsart}
\usepackage{amsmath,amsthm,amssymb,amsfonts}
\usepackage{latexsym}

\DeclareMathOperator{\RE}{Re}

\theoremstyle{definition}
\newtheorem{defn}{Definition}[section]
\theoremstyle{plain}
\newtheorem{rem}[defn]{Remark}
\theoremstyle{definition}

\theoremstyle{plain}

\newtheorem{thm}[defn]{Theorem}

\theoremstyle{plain}
\newtheorem{prop}[defn]{Proposition}
\theoremstyle{plain}
\newtheorem{cor}[defn]{Corollary}

\numberwithin{equation}{section}

\begin{document}
\title[Strong peak points]{Strong peak points and denseness of \\strong peak
functions}
\date{}

\author[H.~J.~Lee]{Han Ju Lee}

\thanks{This work was supported by the Korea Research
Foundation Grant funded by the Korean Government(MOEHRD)
(KRF-2006-352-C00003)}

\baselineskip=.6cm

\begin{abstract} Let $C_b(K)$ be the set of all bounded continuous (real or complex)
functions on a complete metric space $K$ and $A$ a closed subspace
of $C_b(K)$. Using the variational method, it is shown that the set
of all strong peak functions in $A$ is dense if and only if the set
of all strong peak points is a norming subset of $A$. As a corollary
we show that if $X$ is a locally uniformly convex, complex Banach
space, then the set of all strong peak functions in
$\mathcal{A}(B_X)$ is a dense $G_\delta$ subset. Moreover if $X$ is
separable, smooth and locally uniformly convex, then the set of all
norm and numerical strong peak functions in $\mathcal{A}_u(B_X:X)$
is a dense $G_\delta$ subset. In case that a set of uniformly
strongly exposed points of a (real or complex) Banach space $X$ is a
norming subset of $\mathcal{P}({}^n X)$ for some $n\ge 1$, then the
set of all strongly norm attaining elements in $\mathcal{P}({}^n X)$
is dense, in particular, the set of all points at which the norm of
$\mathcal{P}({}^n X)$ is Fr\'echet differentiable is a dense
$G_\delta$ subset.
\end{abstract}

\maketitle

\section{Main Result}

Let $K$ be a complete metric space and $C_b(K)$ the Banach space of
all bounded (real or complex) continuous functions on $K$ with sup
norm $\|f\|= \sup \{ |f(t)| : t \in K\}$. A nonzero function $f\in
A$ is said to be a {\it strong peak function} at $t\in K$ if
whenever there is a sequence $\{t_n\}$ such that $\lim_k |f(t_k)|=
\|f\|$, we get $\lim_k t_k = t$. The corresponding point $t\in K$ is
said to be a {\it strong peak point} for $A$. We denote by $\rho A$
the set of all strong peak points for $A$.

Bishop's theorem \cite{B} says that if $A$ is a (complex) uniform
algebra on a compact metric space $K$, then the set $\rho A$ is a
norming subset for $A$. That is, $\|f\|= \sup_{t\in K} |f(t)|$ for
each $f\in A$. It is observed  \cite{CLS}  that the denseness of the
set of all strong peak functions implies that $\rho A$ is a norming
subset of $A$ and Bishop's theorem is generalized there with
applications to the existence of numerical boundaries.

In this paper, we prove that the converse holds. Precisely if $\rho
A$ is a norming subset for $A$, then the set of all strong peak
functions in $A$ is dense. We use the variational method similar to
\cite{DGZ}.

\begin{thm}\label{thm:main} Let $A$ be a closed subspace of $C_b(K)$, where $K$ is a
complete metric space. The set $\rho A$ is a norming subset of $A$
if and only if the set of all strong peak functions in $A$ is a
dense $G_\delta$ subset of $A$.
\end{thm}
\begin{proof}
Let $d$ be the complete metric on $K$. Fix $f\in A$ and
$\epsilon>0$. For each $n\ge 1$, set
\[ U_n =\{ g\in A: \exists z\in \rho A \mbox{ with }
|(f-g)(z)| > \sup\{ |(f-g)(x)| : d(x,z)> 1/n\} \}.\] Then $U_n$ is
open and dense in $A$. Indeed, fix $h\in A$. Since $\rho A$ is a
norming  subset of $A$, there is a point $w\in \rho A$ such that
\[ |(f-h)(w)| > \|f-h\| -\epsilon/2.\]
Put $g(x) = h(x) + \epsilon\cdot p(x)$, where $p$ is a strong peak
function at $w$ such that $|p(x)|<1/3$ for $\|x-w\|>1/n$, $|p(w)|=1$
and $|f(w)-h(w) - \epsilon p(w)| =|f(w)-h(w)| + \epsilon$. Then
$\|g-h\|\le \epsilon$ and
\begin{align*}
|(f-g)(w)|&=|f(w)-h(w)-\epsilon p(w)| = |f(w)-h(w)|+ \epsilon\\
&> \|f-h\| + \epsilon/2\\
&\ge \sup\{ |(f-h)(x) - \epsilon p(x)|: d(x,w)> 1/n\}.\\
&=\sup\{ \|(f-g)(x)|: d(x, w)>1/n\}.
\end{align*} That is, $g\in U_n$.

By the Baire category theorem there is a $g\in \bigcap U_n$ with
$\|g\|<\epsilon$, and we shall show that $f-g$ is a strong peak
function. Indeed, $g\in U_n$ implies that there is $z_n\in X$ such
that
\[ |(f-g)(z_n)|> \sup\{ |(f-g)(x)|: d(x,z_n)> 1/n\}.\]
Thus $d(z_p, z_n)\le 1/n$ for every $p>n$, and hence $\{z_n\}$
converges to a point $z$, say. Suppose that there is another
sequence $\{x_k\}$ in $B_X$ such that $\{|(f-g)(x_k)|\}$ converges
to $\|f-g\|$. Then for each $n\ge 1$, there is $M_n\ge 1$ such that
for every $m\ge M_n$,
\[ |(f-g)(x_m)|> \sup\{ |(f-g)(x)|: d(x,z_n)> 1/n\}.\] Then $d(x_m, z_n)\le
1/n$ for every $m\ge M_n$. Hence $\{x_m\}$ converges to $z$. This
shows that $f-g$ is a strong peak function at $z$.

The direct argument shows that the set of all  strong peak functions
of $A$ is a $G_\delta$ subset of $A$ (cf. Proposition~2.19 in
\cite{KL}). This proves the necessity.

Concerning the converse, it is shown \cite{CLS} that if the set of
all strong peak functions is dense in $A$, then the set of all
strong peak points is a norming subset of $A$.
\end{proof}

\begin{rem}
Let $Y$ be a Banach space and  $C_b(K:Y)$ the Banach space of all
bounded continuous functions from a complete metric space $K$ into
$Y$ with the sup norm $\|f\|=\sup\{\|f(x)\|:x\in K\}$ for each $f\in
C_b(K:Y)$. Then Theorem~\ref{thm:main} also holds for each closed
subspace $A$ of $C_b(K:Y)$.
\end{rem}

Let $B_X$ be the unit ball of the Banach space $X$. Recall that the
point $x\in B_X$ is said to be a {\it smooth point} if there is a
unique $x^*\in B_{X^*}$ such that $\RE x^*(x)=1$. We denote by ${\rm
sm} (B_X)$ the set of all smooth points of $B_X$. We say that a
Banach space is {\it smooth} if ${\rm sm}(B_X)$ is the unit sphere
$S_X$ of $X$. The following corollary shows that if $\rho A$ is a
norming subset, then the set of smooth points of $B_A$ is dense in
$S_A$.

\begin{cor}\label{cor:main}
Suppose that $K$ is a complete metric space and $A$ is a subspace of
$C_b(K)$. If $\rho A$ is a norming subset of $A$, then the set of
all smooth points of $B_A$ contains a dense $G_\delta$ subset of
$S_A$.
\end{cor}
\begin{proof}
It is shown \cite{CLS} that if $f\in A$ is a strong peak function
and $\|f\|=1$, then $f$ is a smooth point of $B_A$. Then
Theorem~\ref{thm:main} completes the proof.
\end{proof}

\section{Denseness of strongly norm attaining polynomials}

Let $X$ be a Banach space over a scalar field (real or complex)
$\mathbb{F}$ and $X^*$ the dual space of $X$. If $X$ and $Y$ are
Banach spaces, an {\it $N$-homogeneous polynomial} $P$ from $X$ to
$Y$ is a mapping such that there is an $N$-linear (bounded) mapping
$L$ from $X\times \dots \times X$ to $Y$ such that
\[ P(x) = L(x, \dots, x),~~\forall x\in X.\]
${\mathcal P}(^N X:Y)$ denote the Banach space of all
$N$-homogeneous polynomials from $X$ to $Y$, endowed with the
polynomial norm $\|P\|=\sup_{x \in B_X}{\|P(x)\|}$. When
$Y=\mathbb{F}$, ${\mathcal P}(^N X:Y)$ is denoted by
$\mathcal{P}({}^NX)$. We refer to \cite{D} for background on
polynomials. An $N$-homogeneous polynomial $P:X\to Y$ is said to
{\it strongly attain its norm} at $z$ if whenever there is a
sequence $\{x_n\}$ in $B_X$ such that $\lim_n \|P(x_n)\|=\|P\|$, we
get a convergent subsequence of $\{x_n\}$ which converges to
$\lambda z$ for some $\lambda\in S_\mathbb{C}$.

An element $x\in B_X$ is said to be a {\it strongly exposed point}
for $B_X$ if there is a linear functional $f\in B_{X^*}$ such that
$f(x)=1$ and whenever there is a sequence $\{x_n\}$ in $B_X$
satisfying $\lim_n \RE f(x_n)= 1$, we get $\lim_n \|x_n -x\|= 0$. A
set $\{x_{\alpha}\}$ of points on the unit sphere $S_X$ of $X$ is
called {\it uniformly strongly exposed} (u.s.e.), if there are a
function $\delta(\epsilon)$ with $\delta(\epsilon)>0$ for every
$\epsilon
>0$, and a set $\{f_{\alpha}\}$ of elements of norm $1$ in $X^*$ such that
for every $\alpha$, $ f_{\alpha}(x_{\alpha}) = 1$, and for any $x$,
$$\|x\|\le 1 ~\text{and}~ \RE f_{\alpha}(x) \ge 1- \delta(\epsilon)
~\text{imply}~ \|x-x_{\alpha}\| \le \epsilon.$$ In this case we say
that $\{f_{\alpha}\}$ uniformly strongly exposes $\{x_{\alpha}\}$.
Lindenstrauss \cite[Proposition 1]{Li} showed that if $S_X$ is the
closed convex hull of a set of u.s.e. points, then $X$ has property
$A$, that is, for every Banach space $Y$, the set of norm-attaining
elements is dense in $\mathcal{L}(X,Y)$, the Banach space of all
bounded operators of $X$ into $Y$. The following theorem gives
stronger result.

Recall that the norm $\|\ \|$ of a Banach space is said to be {\it
Fr\'echet differentiable } at $x\in X$ if
\[ \lim_{\delta\to 0} \ \sup_{\|y\|=1} \frac{\|x+\delta y\| +
\|x-\delta y\| -2\|x\|}{\delta}=0.\] It is well-known that the set
of Fr\'echet differentiable points in a Banach space is a $G_\delta$
subset \cite[Proposition~4.16]{BL}. It is shown \cite{F} that in a
real Banach space $X$, the norm of $\mathcal{P}({}^n X)$ is
Fr\'echet differentiable at $Q$ if and only if $Q$ strongly attains
its norm.

\begin{thm}
Let $X$, $Y$ be Banach spaces over $\mathbb{F}$ and $k\ge 1$.
Suppose that the u.s.e. points $\{x_\alpha\}$ in $S_X$ is a norming
subset of $\mathcal{P}({}^k X)$. Then the set of strongly norm
attaining elements in $\mathcal{P}({}^k X:Y)$ is dense. In
particular, the set of all points at which the norm of
$\mathcal{P}({}^n X)$ is Fr\'echet differentiable is a dense
$G_\delta$ subset.
\end{thm}

\begin{proof}
Let $\{x_\alpha\}$ be a u.s.e. points and $\{x_\alpha^*\}$ be the
corresponding functional which uniformly strongly exposes
$\{x_\alpha\}$. Let $A=\mathcal{P}({}^k X:Y)$, $f\in A$ and
$\epsilon>0$. Fix $n\ge 1$ and set
\[ U_n =\{ g\in A: \exists z\in \rho A \mbox{ with }
\|(f-g)(z)\| > \sup\{ \|(f-g)(x)\| : \inf_{|\lambda|=1}\|x-\lambda
z\|> 1/n\} \}.\] Then $U_n$ is open and dense in $A$. Indeed, fix
$h\in A$. Since $\{x_\alpha\}$ is a norming  subset of $A$, there is
a point $w\in \{x_\alpha\}$ such that
\[ |(f-h)(w)| > \|f-h\| -\delta(1/n)\epsilon.\]
Put $g(x) = h(x) -\epsilon\cdot
p(x)^k\frac{f(w)-h(w)}{\|f(w)-h(w)\|}$, where $p$ is a strongly
exposed functional at $w$ such that $|p(x)|<1-\delta(1/n)$ for
$\inf_{|\lambda|=1}\|x-\lambda w\|>1/n$, $p(w)=1$. Then $\|g-h\|\le
\epsilon$ and
\begin{align*}
\|(f-g)(w)\|&=\|f(w)-h(w)+\epsilon p(w)^k\frac{f(w)-h(w)}{\|f(w)-h(w)\|}\| = \|f(w)-h(w)\|+ \epsilon\\
&> \|f-h\| + \epsilon(1-\delta(1/n))\\
&\ge \sup\{ \|(f-h)(x) - \epsilon p(x)^k\frac{f(w)-h(w)}{\|f(w)-h(w)\|}\|: \inf_{|\lambda|=1}\|x-\lambda w\|> 1/n\}.\\
&=\sup\{ \|(f-g)(x)|:\inf_{|\lambda|=1}\|x-\lambda w\|> 1/n\}.
\end{align*} That is, $g\in U_n$.

By the Baire category theorem there is a $g\in \bigcap U_n$ with
$\|g\|<\epsilon$, and we shall show that $f-g$ is a strong peak
function. Indeed, $g\in U_n$ implies that there is $z_n\in X$ such
that
\[ \|(f-g)(z_n)\|> \sup\{ \|(f-g)(x)\|: \inf_{|\lambda|=1}\|x-\lambda z_n\|> 1/n\}.\]
Thus $\inf_{|\lambda|=1}\|z_p -\lambda z_n\|\le 1/n$ for every
$p>n$, and $\inf_{|\lambda|=1}|z_p^*(z_n)-\lambda|=
1-|z_n^*(z_p)|\le 1/n$ for every $p>n$. So $\lim_n\inf_{p>n}
|z_n^*(z_p)|=1$. Hence there is a subsequence of $\{z_n\}$ which
converges to $z$, say by \cite[Lemma~6]{A1}. Suppose that there is
another sequence $\{x_k\}$ in $B_X$ such that $\{\|(f-g)(x_k)\|\}$
converges to $\|f-g\|$. Then for each $n\ge 1$, there is $M_n$ such
that $M_n\ge n$ and  for every $m\ge M_n$,
\[ \|(f-g)(x_m)\|> \sup\{ \|(f-g)(x)\|:\inf_{|\lambda|=1} \|x-\lambda z_n\|> 1/n\}.\] Then
$\inf_{|\lambda|=1}\|x_m-\lambda z_n\|\le 1/n$ for every $m\ge M_n$.
So $\inf_{|\lambda|=1}\|x_m-\lambda z\|\le
\inf_{|\lambda|=1}\|x_m-\lambda z_n \| + \|z- z_n \|\le 2/n$ for
every $m\ge M_n$. Hence we get a convergent subsequence of $x_n$ of
which limit is $\lambda z$ for some $\lambda\in S_\mathbb{C}$. This
shows that $f-g$ strongly norm attains at $z$.

It is shown \cite{CLS} that the norm is Fr\'echet differentiable at
 $P$ if and only if whenever there are sequences $\{t_n\}$,
$\{s_n\}$ in $B_X$ and scalars $\alpha$, $\beta$ in $S_\mathbb{F}$
such that $\lim_n \alpha P(t_n)=\lim_n \beta P(s_n) = \|P\|$, we get
\begin{equation}\label{eq1}\lim_n\ \sup_{\|Q\|=1}(\alpha Q(t_n)-\beta Q(s_n))=0.\end{equation}
We have only to show that every nonzero element $P$ in $A$ which
strongly attains its norm satisfies the condition~(\ref{eq1}).
Suppose that $P$ strongly attains its norm at $z$ and $P\neq 0$.

For each $Q\in A$, there is a $k$-linear form $L$ such that $Q(x) =
L(x,\dots,x)$ for each $x\in X$. The polarization identity \cite{D}
shows that $\|Q\|\le \|L\|\le (k^k/k!)\|Q\|$. Then for each $x,y\in
B_X$, $\|Q(x)- Q(y)\|\le n\|L\|\|x-y\|$ and
\[\|Q(x)- Q(y)\|\le \frac{k^{k+1}}{k!} \|Q\|\|x-y\|.\]

Suppose that there are sequences $\{t_n\}$, $\{s_n\}$ in $B_X$ and
scalars $\alpha$, $\beta$ in $S_\mathbb{F}$ such that $\lim_n \alpha
P(t_n)=\lim_n \beta P(s_n) = \|P\|$, then for any subsequences
$\{s_n'\}$ of $\{s_n\}$ and $\{t_n'\}$ of $\{t_n\}$, there are
convergent further subsequences $\{t''_n\}$ of $\{t_n'\}$ and
$\{s_n''\}$ of $\{s_n'\}$ and scalars $\alpha''$ and $\beta''$ in
$S_\mathbb{F}$ such that $\lim_n t''_n = \alpha'' z$ and $\lim_n
s_n'' = \beta'' z$. Then $\alpha P(\alpha''z) = \beta P(\beta''z) =
1$. So $\alpha(\alpha'')^k = \beta(\beta'')^k$.

Then we get
\begin{align*}\overline{\lim_n}\ \sup_{\|Q\|=1} (\alpha Q(t''_n) - \beta
Q(s''_n)) &\le \overline{\lim_n}\ \sup_{\|Q\|=1} (\alpha Q(t''_n) -
\alpha
Q(\alpha''z)) -(\beta Q(\beta'' z) - \beta Q(s''_n))\\
&\le \overline{\lim_n}  \frac{k^{k+1}}{k!} (\|t''_n - \alpha''z\| +
\| \beta'' z-  s''_n\|)=0.
\end{align*}
This implies that $\lim_n \sup_{\|Q\|=1}(\alpha Q(t_n)- \beta
Q(s_n))=0$. Therefore the norm is Fr\'echet differentiable at $P$.
This completes the proof.
\end{proof}

\begin{rem}
Suppose that the $B_X$ is the closed convex hull of a set of  u.s.e
points, then the set of u.s.e. points is a norming subset of
$X^*=\mathcal{P}({}^1 X)$. Hence the elements in $X^*$ at which the
norm of $X^*$ is Fr\'echet differentiable is a dense $G_\delta$
subset.
\end{rem}

\section{Denseness of strong peak holomorphic functions}

Let $X$ be a complex Banach space and $B_X$ the unit ball of $X$. We
consider two Banach algebras as subspaces of complex $C_b(B_X)$:
\begin{align*} \mathcal{A}_b(B_X) &=\left\{ f\in C_b(B_X): f \mbox{
is holomorphic on the interior of } B_X  \right\}\\
\mathcal{A}_u(B_X) &= \left\{ f\in \mathcal{A}_b(B_X) : f \mbox{ is
uniformly continuous on } B_X \right\}\end{align*} We shall denote
by $\mathcal{A}(B_X)$ either $\mathcal{A}_b(B_X)$ or
$\mathcal{A}_u(B_X)$.

Using the Bourgain-Stegall variational method \cite{B,S}, it is
shown \cite{CLS} that if $X$ is a complex Banach space with the
Radon-Nikod\'ym property,  then the set of all strong peak functions
in $\mathcal{A}(B_X)$ is dense. In case that $X$ is locally
uniformly convex, it is shown \cite{CLS} that the set of norm
attaining elements is dense in $\mathcal{A}(B_X)$. That is, the set
consisting of $f\in \mathcal A(B_X)$ with $\|f\|= |f(x)|$ for some
$x\in B_X$ is dense in $\mathcal{A}(B_X)$.

The following corollary gives a stronger result. Notice that if $X$
is locally uniformly convex, then every point of $S_X$ is the strong
peak point for $A(B_X)$ \cite{G1}. Theorem~\ref{thm:main} and
Corollary~\ref{cor:main} implies the following.
\begin{cor}\label{cor:localconvex}
Suppose that $X$ is a  locally uniformly convex, complex Banach
space. Then the set of all strong peak functions in
$\mathcal{A}(B_X)$ is a dense $G_\delta$ subset of
$\mathcal{A}(B_X)$. In particular, the set of all smooth points of
$B_{\mathcal{A}(B_X)}$ contains a dense $G_\delta$ subset of
$S_{\mathcal{A}(B_X)}$.
\end{cor}

Let $X$ be a complex Banach space. Let $\mathcal{A}_b(B_X:X)$ to be
the Banach space consisting of $X$-valued bounded continuous
functions $f$ from $B_X$ to $X$, which is also holomorphic on the
interior of $B_X$ with the sup norm
\[ \|f\|= \sup\{ \|f(x)\|:x\in B_X\}.\] The space $\mathcal
A_u(B_X:X)$ is the subspace of $\mathcal{A}_b(B_X:X)$ consisting of
all uniformly continuous functions on $B_X$. We denote by $\mathcal
A(B_X:X)$ either $\mathcal{A}_b(B_X:X)$ or $\mathcal{A}_u(B_X:X)$.

Let $\Pi(X)=\{(x,x^*)\in B_X\times B_{X^*}: x^*(x)=1\}$ be the
topological subspace of $B_X\times B_{X^*}$, where $B_X$ and
$B_{X^*}$ is equipped with norm and weak-$*$ compact topology
respectively. For $f\in \mathcal{A}(B_X:X)$, the {\it numerical
radius} $v(f)$ of $f$ is defined by $v(f) = \sup\{ |x^*f(x)|:(x,
x^*)\in \Pi(X)\}$.

The {\it numerical  strong peak function} is introduced in \cite{KL}
and the denseness of holomorphic numerical strong peak functions in
$\mathcal A(B_X:X)$ is studied. The function $f\in
\mathcal{A}(B_X:X)$ is said to be a numerical strong peak function
if there is $(x,x^*)$ such that $\lim_n |x^*_nf(x_n)|=v(f)$ for some
$\{ ( x_n, x_n^*)\}_n$ in $\Pi(X)$ implies that $(x_n, x_n^*)$
converges to $(x, x^*)$ in $\Pi(X)$.  The function $f\in
\mathcal{A}(B_X:X)$ is said to be {\it numerical radius attaining}
if there is $(x, x^*)$ in $\Pi(X)$ such that $v(f) = |x^*f(x)|$.

\begin{prop}\label{prop:numerical}
Suppose that the set $\Pi(X)$ is complete metrizable and the set
$\Gamma= \{ (x, x^*)\in \Pi(X) : x\in \rho \mathcal{A}(B_X)\cap {\rm
sm}(B_X)\}$ is a numerical boundary. That is, $v(f)
=\sup_{(x,x^*)\in \Gamma} |x^*f(x)|$  for each $f\in
\mathcal{A}(B_X:X)$.  Then the set of all numerical strong peak
functions in $\mathcal{A}(B_X:X)$ is a dense $G_\delta$ subset of
$\mathcal{A}(B_X:X)$.
\end{prop}
\begin{proof}
Let $A=\mathcal{A}(B_X:X)$ and let $d$ be a complete metric on
$\Pi(X)$. In \cite{KL}, it is shown that if $\Pi(X)$ is complete
metrizable, then the set of all numerical peak functions in $A$ is a
$G_\delta$ subset of $A$. We need prove the denseness. Let $f\in A$
and $\epsilon>0$. Fix $n\ge 1$ and set
\begin{align*} U_n =\{ g\in A&: \exists (z,z^*)\in
\Gamma \mbox{ with }\\
&|z^*(f-g)(z)| > \sup\{ |x^*(f-g)(x)| : d((x,x^*),(z,z^*))> 1/n\}
\}.\end{align*}Then $U_n$ is open and dense in $A$. Indeed, fix
$h\in A$. Since $\Gamma$ is a numerical boundary for $A$, there is a
point $(w, w^*)\in \Gamma$ such that
\[ |w^*(f-h)(w)| > v(f-h) -\epsilon/2.\]
Notice that $d((x,x^*), (w,w^*))>1/n$ implies that there is
$\delta_n>0$ such that $\|x-w\|>\delta_n$. Choose a peak function
$p\in \mathcal{A}(B_X)$ such that $\|p\|=1=|p(w)|$ and $|p(x)|<1/3$
for $\|x-w\|>\delta_n$ and $|w^*(f-h)(w) - \epsilon p(w)|
=|w^*f(w)-w^*h(w)| + \epsilon|p(w)|= |w^*f(w)-w^*h(w)| + \epsilon$.

Put $g(x) = h(x) + \epsilon\cdot p(x)w$, where $p$ is a peak
function at $w$ such that $|p(x)|<1/3$ for $\|x-w\|>\delta_n$,
\begin{align*}
|w^*(f-g)(w)|&=|w^*f(w)-w^*h(w)-\epsilon p(w)|= |w^*f(w)-w^*h(w)|+ \epsilon\\
&>v(f-h) + \epsilon/2\\
&\ge \sup\{ |x^*(f-h)(x) - \epsilon p(x)x^*(w)|: d((x,x^*), (w,w^*))> 1/n\}.\\
&=\sup\{ |x^*(f-g)(x)|: d((x,x^*), (w,w^*))>1/n\}.
\end{align*} That is, $g\in U_n$.

By the Baire category theorem there is a $g\in \bigcap U_n$ with
$\|g\|<\epsilon$, and we shall show that $f-g$ is a strong peak
function. Indeed, $g\in U_n$ implies that there is $(z_n, z_n^*)\in
\Gamma$ such that
\[ |z_n^*(f-g)(z_n)|> \sup\{ |x^*(f-g)(x)|: d((x,x^*), (z_n,z_n^*))> 1/n\}.\]
Thus $d((z_p,z_p^*), (z_n,z_n^*))\le 1/n$ for every $p>n$, and hence
$\{(z_n, z_n^*)\}$ converges to a point $(z, z^*)$, say. Suppose
that there is another sequence $\{(x_k, x_k^*)\}$ in $\Pi(X)$ such
that $\{|x_k^*(f-g)(x_k)|\}$ converges to $v(f-g)$. Then for each
$n\ge 1$, there is $M_n\ge 1$ such that for every $m\ge M_n$,
\[ |x_m^*(f-g)(x_m)|> \sup\{ |x^*(f-g)(x)|: d((x,x^*), (z_n, z_n^*))> 1/n\}.\]
Then $d((x_m, x_m^*), (z_n, z_n^*))\le 1/n$ for every $m\ge M_n$.
Hence $\{(x_m, x_m^*)\}$ converges to $(z, z^*)$. This shows that
$f-g$ is a numerical strong peak function at $(z, z^*)$.
\end{proof}

It is shown \cite{KL} that the set of all numerical radius attaining
elements is dense in $\mathcal{A}(B_X:X)$ if $X$ is locally
uniformly convex. In addition, if $X$ is separable and smooth, we
get the stronger result on $\mathcal{A}_u(B_X:X)$.

\begin{cor}
Suppose that $X$ is separable, smooth and locally uniformly convex.
Then the set of norm and numerical strong peak functions is a dense
$G_\delta$ subset of $\mathcal{A}_u(B_X:X)$.
\end{cor}
\begin{proof}
It is shown \cite{KL} that if $X$ is separable, $\Pi(X)$ is complete
metrizable. In view of \cite[Theorem~2.5]{Pal}, $\Gamma$ is a
numerical boundary for $\mathcal{A}_u(B_X:X)$. Hence
Proposition~\ref{prop:numerical} shows that the set of all numerical
strong peak functions is dense in $\mathcal{A}_u(B_X:X)$. Finally,
Corollary~\ref{cor:localconvex} implies that the set of all norm and
numerical peak functions is a dense $G_\delta$ subset of
$\mathcal{A}_u(B_X:X)$.
\end{proof}

It is shown \cite{CHL} that the set of all strong peak points for
$\mathcal{A}(B_X)$ is dense in $S_X$ if $X$ is an order continuous
locally uniformly $c$-convex sequence space. (For the definition see
\cite{CHL}).

\begin{cor}
Let $X$ be an order continuous locally uniformly $c$-convex, smooth
Banach sequence space. Then the set of all norm and numerical strong
peak functions in $\mathcal{A}_u(B_X:X)$ is a dense $G_\delta$
subset of $\mathcal{A}_u(B_X:X)$.
\end{cor}
\begin{proof}
Notice that $X$ is separable since $X$ is order continuous. Hence
the set of all smooth points of $B_X$ is dense in $S_X$ by the Mazur
theorem and $\Pi(X)$ is complete metrizable \cite{KL}. In view of
\cite[Theorem~2.5]{Pal}, $\Gamma$ is a numerical boundary for
$\mathcal{A}_u(B_X:X)$. Hence Proposition~\ref{prop:numerical} shows
that the set of all numerical strong peak functions is a dense
$G_\delta$ subset of $\mathcal{A}_u(B_X:X)$. Theorem~\ref{thm:main}
also shows that the set of all strong peak functions is a dense
$G_\delta$ subset of $\mathcal{A}_u(B_X:X)$. This completes the
proof.
\end{proof}

\bibliographystyle{amsplain}

{\noindent Mathematics Department\\
 University of Missouri-Columbia\\
e-mail:hahnju@postech.ac.kr}

\end{document}